\documentclass[%
 aip,
 amsmath,amssymb,
 preprint,%
]{revtex4-2}

\usepackage{graphicx}
\usepackage{dcolumn}
\usepackage{bm}

\usepackage[utf8]{inputenc}
\usepackage[T1]{fontenc}
\usepackage{mathptmx} 

\usepackage{subfig}
\usepackage{mwe}
\usepackage[percent]{overpic}
\usepackage{amssymb,amsmath,epsfig,mathrsfs,algorithm,algorithmic,rotating}
\usepackage{mathtools,graphics,psfrag,graphicx,color,stmaryrd}
\usepackage[table]{xcolor}
\usepackage{tikz}
\usetikzlibrary{arrows,shapes,backgrounds}
\usepackage{booktabs}
\usepackage{pgfplots}
\pgfplotsset{compat=newest}
\usepackage{slashed}
\usepackage{pgfplotstable}
\usepackage{array}
\usepackage{colortbl}
\usepackage{booktabs}
\usepackage{eurosym}
\usepackage{amsmath}
\usepackage{pgfplotstable}
\usepackage{pgfplots}
\pgfplotsset{compat=newest}
\usepackage{verbatim}
\usepackage[normalem]{ulem} 
\usepackage{fancyhdr}
\usepackage{hyperref}
\hypersetup{
pagebackref=true, 
hyperindex=true, 
colorlinks=true, 
breaklinks=true, 
urlcolor=blue, 
linkcolor=blue, 
citecolor=blue, 
filecolor=blue,
pdftitle={Composition of one-step methods}, 
pdfauthor={Aymen Laadhari}, 
pdfsubject={Composition of one-step methods} 
}

\newcommand{\diam}{\operatorname*{diam}}

\newcommand\ba{\boldsymbol{y}}
\newcommand{\iu}{{i\mkern1mu}}
\newcommand\bx{\boldsymbol{x}}
\newcommand\bu{\boldsymbol{u}}

\newcommand{\Frac}[2] {\displaystyle \frac{#1}{#2}}

\newcommand{\Th}{\mathscr{T}_{h}}
\newcommand{\T}{\mathscr{T}_{h}}

\def\bnu{\boldsymbol{\nu}}

\def\G{\Gamma}

\def\bI{\boldsymbol{I}}

\def\by{\boldsymbol{y}}

\def\G{\Gamma}

\def\bI{\boldsymbol{I}}

\def\by{\boldsymbol{y}}

\newtheorem{theorem}{Theorem}

\usepackage{marginnote}

\newcommand{\cblue}[1]{\textcolor{black}{#1}}
\newcommand{\cred}[1]{\textcolor{black}{#1}}

\colorlet{cgray}{gray!20!white}
\colorlet{cblue}{blue!20!white}
\colorlet{RoyalBlue}{blue!20!white}
\colorlet{RoyalGreen}{green!20!white}

\newcommand{\bnabla}    { \boldsymbol{\nabla} }

\renewcommand{\leq}    {\leqslant}
\renewcommand{\geq}    {\geqslant}

\everymath{\displaystyle}

\allowdisplaybreaks

\begin{document}

\title{\cblue{Numerical Approach Based on the Composition of One-Step Time-Integration Schemes For Highly Deformable Interfaces}}

\author{Aymen Laadhari}%
\author{Ahmad Deeb} 
\affiliation{Department of Mathematics, College of Arts and Sciences, Khalifa University of Science and Technology, Abu Dhabi, United Arab Emirates.
}

\date{December 2, 2022} 

\begin{abstract}
In this work, we propose a numerical approach for simulations of large deformations of interfaces in a level set framework. To obtain a fast and viable numerical solution in both time and space, temporal discretization is based on the composition of one-step methods exhibiting higher orders and stability, especially in the case of stiff problems with strongly oscillatory solutions. Numerical results are provided in the case of ordinary and partial differential equations to show the main features and demonstrate the performance of the method. Convergence properties and efficiency in terms of computational cost are also investigated.
\end{abstract}

\maketitle

\section{Introduction}

We are interested in numerical solutions of strongly coupled systems of PDEs involving highly nonlinear membrane forces, see e.g.~\cite{LAADHARI201835,LAADHARI2017}.
For such problems with high nonlinearities, fine meshes and spatial discretizations by high-order finite elements are required~\cite{Doyeux2013251}. Therefore, numerical time integrations with higher orders are also necessary to allow viable numerical approximations.
In this preliminary work, we present the formalism in the case of the level set problem, where we seek to produce solutions with a higher-order time approximations.

Different temporal integration techniques with high order approximations exist to study dynamical systems.
Popular numerical methods are those associated with discrete flow, where solutions are sought over discrete times using a numerical flow composition \cite{book:hairer}. Fixed or adapted time steps were used.
In the existing literature, several studies on continuous flows have investigated the Borel-Pad\'e-Laplace integrator for stiff and non-stiff problems. See, for example, \cite{book:hairer2,ahmad_comp_bpl_sfg_2015}. Indeed, an integrator calculates an approximation over a continuous time interval and has the ability to increase the order of approximation as much as necessary by changing a single parameter in the integrator \cite{DEEB_2022_bpl}.

To allow computational savings and improve precision while maintaining numerical accuracy, adaptive time-stepping strategies are commonly employed where the level of adaptability may depend on the rate of change between consecutive solutions~\cite{GUO2022107943} or on some error estimates \cite{book:hairer,book:hairer2}. However, the computational efficiency depends on the order of the scheme and how the error estimate is computed.
Hence the need for high-order discrete flow to improve numerical efficiency, especially for stiff problems. For a given numerical approximation, the convergence depends on the error in time and in space. Although the adaptive meshing techniques allow a better precision of the numerical solution, the global error estimate will be limited by the order of the numerical scheme in time. Thus, it is important to increase the order of approximation in time to allow high-order of precision.

One can distinguish one-step and multi-step methods based on the number of previous numerical solutions used in the temporal discretization of differential equation. 
For a detailed discussion of the classical theory of multi-step methods as well as the properties of order, convergence, stability and symmetry, we refer the reader to \cite{book:butcher,book:hairer2,hairer2002geometric}.
Among the most used schemes for stiff problems with  oscillatory terms, we can mention the implicit linear multi-step methods based on the backward difference formulas (BDF) presenting a stability limited to lower orders \cite{iserles_2008,curtiss_1952,butcher_1963,Butcher_1996}. Indeed, only first-order (implicit Euler) and second-order schemes are stable. The BDF schemes of order 3 to 6 exhibit a weaker stability property, called A-stability, and formulas of order greater than 6 are unstable. A generalization of BDF using a second derivative allows to obtain implicitly stable schemes up to the order 10, see e.g. \cite{book:hairer2}.
Runge-Kutta methods, referred to as RK, are one-step schemes \cite{Runge_1895}. 
\cblue{To cope with the lack of A stability of explicit RK methods, which are not suitable for solving stiff equations, implicit and stable RK methods are also developed. We refer the reader to \cite{book:butcher,book:hairer2} for a note on A-stability and the Gauss, Radau IIA and Lobatto IIIA, IIIB and IIIC schemes.}
Composition methods have emerged in the literature, with several variants, as a powerful numerical tool based on the composition of lower-order basic time-integration methods with the aim of increasing the order of approximation.

The aim of this work is to provide a higher-order approximation of the level set solution. Various higher-order finite element approximations are used and allow to increase the spatial precision. 
Since the error in time takes over in case of low-order schemes, composition methods are applied on simple and low-order discrete flows to increase the order of approximation in time, while keeping a reasonable computational cost.
Numerical examples are presented to validate and show the main characteristics of the method in the cases of the Lotka-Volterra differential system and the level set advection problem.
This paper is organized as follows.
Section 1 presents the preliminary concepts and mathematical setting.
In Section~2, we present the composition methods for the implicit  first-order accurate backward Euler and  second-order accurate Heun's method.
In Section~3, we provide an assessment of the composition technique through numerical examples.

%
\vspace{-0.1cm}
\section{{Mathematical formulation}}\label{sec:model}
\subsection{Level set problem}
Let $T>0$ and $d=2$ be respectively the simulation period and the spatial dimension.
For any time $t\in(0,T)$, the interface is denoted by $\Gamma(t)$.
Its deformations are described implicitly as the iso-value zero of a  level set function $\varphi$.
The interface is assumed smooth enough and is enclosed in a larger domain $\Lambda$ so that
$$\G(t)= \big\lbrace(t,\bx) \in (0,T)\times\Lambda: \varphi(t,\bx) = 0 \big\rbrace.$$
Given a velocity field $\bu$, the evolution of the level set function satisfies the advection equation:
\begin{equation}
\partial_t\varphi + \bu\cdot\bnabla\varphi = 0, \quad \mbox{ in } (0,T)\times\Lambda.
\label{eq-transport-1}
\end{equation}
Appropriate initial and boundary conditions are considered:
$\varphi = \varphi_b$ on $(0,T)\times \Sigma_{-}$ and $\varphi(0) = \varphi_0$ in $\Lambda$,
where $\Sigma_{-} = \{\bx \in \partial\Lambda: \bu\cdot\bnu(\bx) < 0\}$
is the inflow boundary and $\bnu$ denotes the unit normal exterior vector to $\Lambda$.
For an  initial internal domain $\Omega_0$,
$\varphi$ is initialized as a signed distance function to $\G(t=0)$ so that
$|\bnabla\varphi_0|=1$:
$$ \varphi_0(\bx)=  \left\{
   \begin{array}{ll}\displaystyle
       \inf\{ |\by-\bx| ; \by\in\partial\Omega_0 \}, & \qquad  \text{ if } \bx \notin \Omega_0 , \\
       -\inf\{|\by-\bx| ; \by\in\partial\Omega_0\},& \qquad  \text{ if not.}
 \end{array} \right.
$$
However, the signed distance property is not preserved when solving~\eqref{eq-transport-1},
thus leading to numerical instabilities if the gradient of the level set function becomes either very small, or very large, in particular in the vicinity of the interface or on 
$\partial\Lambda$~\cite{OSHER2001,these}.
To restore the signed distance property, an auxiliary redistancing problem, is commonly solved 
       while maintaining the zero-level set position to avoid the loss of mass characterizing Eulerian methods.
Higher-order methods and the use of fine meshes in the vicinity of the interface can help enforcing the local mass conservation, see for example~\cite{Doyeux2013251,LAADHARI2018376,Laadhari20171047}.

\subsection{Time integration by composition method}

\cblue{Let $y(t)$ be the solution of the following initial value problem:}
\begin{equation}
\label{ode1}
 y^{\prime} =  f\,\left(t,y\,(t)\right),\,\,\text{ for all}\,\,  t\in \mathbb{R},  \text{ with the initial condition } y(0) = y_0.
\end{equation}
For $t>0$, the exact flow of \eqref{ode1} is defined by the map $\mathcal{Y}_t$ as follows:
\begin{equation}
\mathcal{Y}_t:
 \begin{array}{ccl}
 \mathbb{R} &\rightarrow &\mathbb{R}\\
 y_0 &\mapsto & {\mathcal{Y}}_t(y_0) = y(t).
 \end{array}
\end{equation}
For any given numerical scheme of order $p$, we can associate a numerical flow, denoted \cblue{$\varPhi_{\Delta t}^{[p]}$}, such that a sequence of numerical values $y_n$ is constructed on a discrete set of points $ t^n = n\Delta t$ as the approximations of the solution at time $t^n$:
$y_{n+1} =  \varPhi^{[p]}_{\Delta t}\big(y_n \big)$.
We say that the numerical flow is of order $p$ if the following equality holds:
\begin{equation*}
 \mathcal{Y}_{\Delta t}(y_0) - \varPhi^{[p]}_{\Delta t}(y_0)  = O(\Delta t^{p+1}).
\end{equation*}
The superscript $[p]$ in  $\varPhi^{[p]}_{\Delta t}$ refers to the order $p$ of the numerical flow.
As stated in \cite[Theorem  4.1, page 43]{hairer2002geometric}, we can increase the order of the approximation using a $s$ times composition of the same numerical flow $\varPhi^{[p]}_{\Delta t}$.
The resulting discrete flow
$\varPhi^{[p+1]}_{\Delta t}  = \varPhi^{[p]}_{a_1 \Delta t}\circ \varPhi^{[p]}_{a_2 \Delta t} \circ \ldots\circ \varPhi^{[p]}_{a_s \Delta t}$
is at least or order $p+1$ if the conditions hold:
\begin{equation}
\label{comp_flow_cond}
 a_1 + a_2 +\dots +a_s = 1  \qquad\text{ and } \qquad a_1^{p+1} +a_2^{p+1}+ \dots a_s^{p+1}  =0.
\end{equation}
There are no real solutions for system \eqref{comp_flow_cond} when $p$ is odd. However, if $p$ is even, the smallest composition that increase the order by $1$ is when $s=3$. This composition is called the { the triple jump}, see \cite[page 44, section II.4]{hairer2002geometric}. One can also continue to increase the order of construction of the discrete flow by continuing to compose discrete flows, see \cite{murua_1999}. However, a problem for compositions with real coefficients comes from negative coefficients which result in a negative time step. Castella et al. consider solutions for the system with complex coefficients in the case of parabolic equations \cite{castella_2009}.
In the following, we will use the composition technique in the case $s=2$.

Let us consider any discrete flow $\varPhi^{{[p]}}_{\Delta t}$ of order $p$ and introduce the following composition:
$ 
 \varPhi^{[p+1]}_{\Delta t} = \varPhi^{{[p]}}_{a_2^{[p]}\Delta t}\circ \varPhi^{[p]}_{a_1^{[p]}\Delta t}$.
We recall a useful theorem which states that (see e.g. \cite{casas_2021_complex}):
\\
\begin{theorem}
{\textit{Let us consider a differential system of type \eqref{ode1} and any numerical flow $\varPhi^{[p]}_{\Delta t}$ of order $p$. If
 \begin{equation}
\label{coeff_ai_composition}
\overline{a_2^{[p]}} = a_1^{[p]} = \frac{1}{2} + \frac{\iu}{2} \Frac{\sin\left(\Frac{2l+1}{p+1}\pi\right)}{1+\cos\left(\Frac{2l+1}{p+1}\pi\right)},
\qquad
\cred{
\text{ with } \quad 
\left\lbrace
\begin{array}{ll}
  -\frac{p}{2} \leq l \leq\frac{p}{2} -1 & \text{if $p$ is even},\\
 -\frac{p+1}{2} \leq l \leq\frac{p-1}{2} & \text{if $p$ is odd.}
\end{array}
\right.}
 \end{equation}
then the approximation defined by the real part of the output, denoted by $\operatorname{Re}(.)$ and obtained by the composition of flow, is an approximation of the solution of order $p+1$, i.e.:
\begin{equation}
\label{compo_2_BE_cond}
 \mathcal{Y}_{\Delta t}(y_0) - \operatorname{Re}\left( \varPhi^{{[p]}}_{a_2^{[p]}\Delta t}\circ \varPhi^{[p]}_{a_1^{[p]}\Delta t} (y_0)\right) = O\left(\Delta t^{p+2}\right).
\end{equation}
}}
\end{theorem}

First, we consider the numerical flow $\varPhi^{\text{BE1}}$ of the backward Euler scheme which is of order one.
We define the two time composition as follows:
\begin{equation}
\label{compo_2_BE}
 \varPhi^{\text{BE2}}_{\Delta t} :=\varPhi^{\text{BE1}}_{a_2^{[1]}\Delta t} \circ\varPhi^{\text{BE1}}_{a_1^{[1]}\Delta t}.
\end{equation}
Here, the coefficients $a_i^{[1]}$ are given by \cred{formula \eqref{coeff_ai_composition} for $l=0$: } $a_1^{[1]} = \overline{a_2^{[1]}} = \frac{1}{2} + \frac{\iu}{2}$.
Hence, the approximation given by the real part of the output obtained by the composition \eqref{compo_2_BE} is an approximation of the solution of order $2$, i.e.:
\begin{equation*}
\label{compo_2_BE_error}
 \mathcal{Y}_{\Delta t}(y_0) - \operatorname{Re}\left( \varPhi^{\text{BE2}}_{\Delta t}  (y_0)\right) = O\left(\Delta t^3\right).
\end{equation*}
Thus, the double composition of the discrete scheme BE1 of order 1 gives a scheme of order 2.

Second, we consider a numerical flow of order 2 associated to the Heun's method $\varPhi^{HM1}_{\Delta t}$ (second-order RK):
\begin{equation*}
  \mathcal{Y}_{\Delta t}(y_0) - \varPhi^{\text{HM1}}_{\Delta t}  (y_0) = O\left(\Delta t^3\right),
\end{equation*}
thus the composition method:
$\displaystyle
 \varPhi^{\text{HM2}}_{\Delta t} :=\varPhi^{\text{HM1}}_{a_2^{[2]}\Delta t} \circ\varPhi^{\text{HM1}}_{a_1^{[2]}\Delta t},
$ with $6a_1^{[2]} = 6\overline{a_2^{[2]}} = 3 + \iu\sqrt{3}$ 
leads to a numerical flow of order $3$ and we have:
\begin{equation*}
\label{compo_2_HM_error}
 \mathcal{Y}_{\Delta t}(y_0) - \operatorname{Re}\left( \varPhi^{\text{HM2}}_{\Delta t}  (y_0)\right) = O\left(\Delta t^4\right).
\end{equation*}

We can continue increasing the order of the approximation by constructing new numerical flows using the above composition technique. This was thoroughly analysed in \cite{casas_2021_complex} and presented with a comprehensive study.
In addition, it is demonstrated that the numerical flow $\frac{1}{2} \operatorname{Re}\left( \varPhi^{[p]}_{a_1^{[p]} \Delta t}\circ \varPhi^{[p]}_{a_2^{[p]} \Delta t} +\varPhi^{[p]}_{a_2^{[p]} \Delta t}\circ \varPhi^{[p]}_{a_1^{[p]} \Delta t} \right)$ is of order $p+2$.
When composing two discrete flows with permuted coefficients, the imaginary parts are conjugated and are therefore canceled by the summation, resulting in only a real part.

\cred{Remark that the A-stability of the resulting composition scheme is not studied here but we refer the interested reader to the reference \cite{math10224327}.
%
Roughly speaking, the stability function $S_r(z)$ of the resulting composition scheme is the product of two functions: $S_r(z)= S(a_2z)\,S(a_1z)$, where $S(z)$ represents the stability function of the initial diagram. Therefore, the region of stability, defined by $\left\lbrace z \in \mathbb{C},\,\mbox{ such that }\, |S_r(z)|\leq 1 \right\rbrace$, can be represented as the following intersection: $\left\lbrace z \in \mathbb{C}:\, |S(a_2z)|\leq 1 \right\rbrace \cap \left\lbrace z \in \mathbb{C}: \, | S(a_1z)|\leq 1 \right\rbrace$.
Note that each set in this intersection has a wider region than the original scheme with a slight translation in the complex plane.
}

\section{Numerical approximation}\label{sec:method}
\paragraph{\bf Space discretization by finite elements.} 

Let $\Th$ be  a partition of $\Lambda$ consisting of geometrically conforming
open simplicial elements $K$ (triangles for
$d=2$), such that $\displaystyle\overline{\Lambda}=\displaystyle\cup_{K\in\Th}$. 
Lagrange finite element polynomials are considered for the space discretization of 
the level set function on a mesh of size $h=\max_{K\in\Th}\diam(K)$.  Let us denote by $\varphi_{h}$  an
approximation of $\varphi$,
while
$\mathbb{X}_h = \left\lbrace \psi_h\in C^0(\overline{\Lambda}): \psi_{h|K}\in\mathbb{P}^\kappa, \forall
  K\in\Th \right\rbrace$,  with $\kappa\geq 1$
represents  the finite dimensional space of admissible level set.
The Galerkin scheme associated with~\eqref{eq-transport-1} consists 
in finding  $\varphi_h\in \mathbb{X}_h$ such that
$
\displaystyle\varphi_h(t,\bx) = \sum\limits_i y_i(t)\,\psi_{h,i}(\bx)
$, for a certain finite element basis.
\\
A  semi discrete equation is then obtained by multiplying by a test function and integrating aver $\Lambda$.
We employ the Stabilization Upwind Petrov Galerkin, referred to as SUPG~\cite{Loch:229393}, method where the SUPG test function is defined as
$v_{h,i}= \psi_{h,i} + \tau_K \,\bu\cdot \bnabla \psi_{h,i}$ 
in order to add a diffusion in the  streamline direction. 
We set a streamline diffusion parameter proportional to the local mesh size so that
\[\tau_K = \frac{C \,h_K }{ \max\left\{|\bu|_{0,\infty,K}, \text{tol}/h_K\right\}},
\]
where $C$ represents a scaling constant and the term $\text{tol}/h_K$ avoids the division by zero.
The semi-discrete system arising from the discretization with finite elements writes:
\begin{equation}
\label{semi_discrete}
 {\bf{M}}(t)\, \ba^{\prime} +  {\bf{K}}(t)\,\ba ={ \bf{0}}, \quad\text{with } \qquad 
{\ba} = (y_1,\ldots,y_i,\ldots)^T.
\end{equation}
Here, $ {\bf{M}}(t)$ is the stabilized mass matrix, while $ {\bf{K}}(t)$ represents the transport matrix, both time-dependent;
$\ba$ is the unknown degrees of freedom for the level set field.
In the sequel, we omit the subscript $h$ referring to the spatial approximation to alleviate the notations.

\vspace{0.20cm}
\paragraph{\textbf{Time advancing scheme.}}\label{ssec:t} 

\begin{figure}[!b]
 \centerline{
     \includegraphics[height=5.8cm]{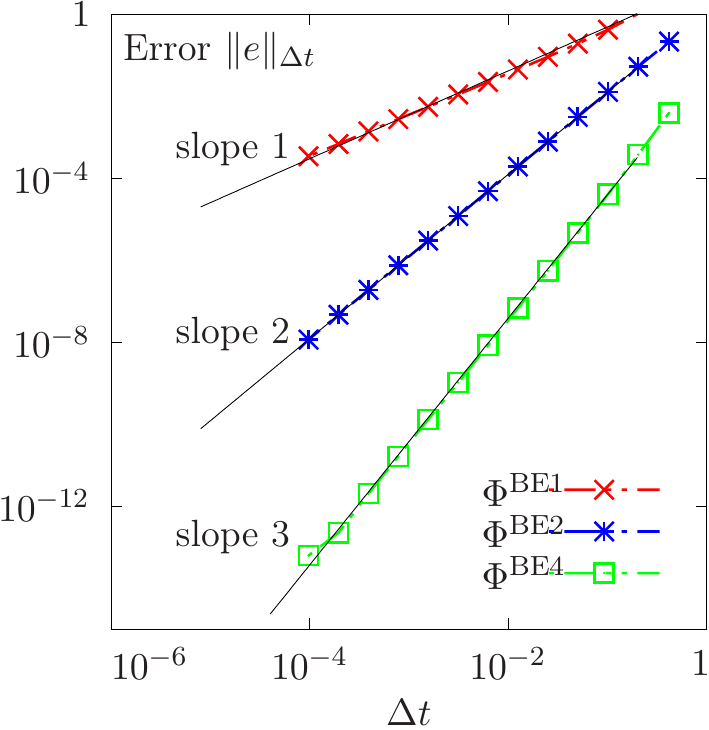}
    \hspace{1.25cm}
     \includegraphics[height=5.8cm]{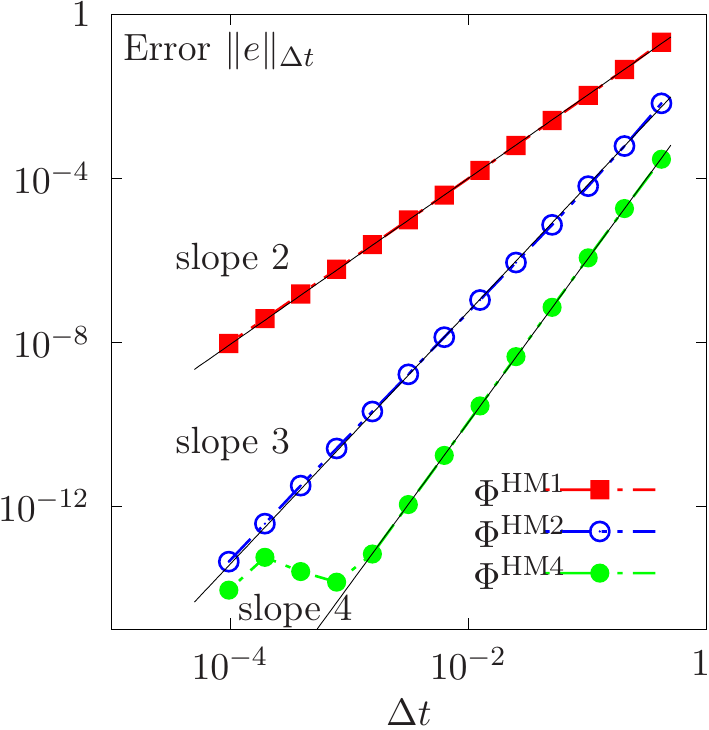}}
 \caption{Error versus the time step size for the composition of backward Euler scheme (left) and Heun's scheme (right).}
 \label{fig:error-t}
\end{figure}

Let us divide  $[0,T]$ into $N$ subintervals $\big[t^n,t^{n+1}\big]$ of constant time step $\Delta t$ with $n=0,..., N-1$.
For $n>0$, the unknown $\varphi_{h}^n$ approximates $\varphi$ at $t^n$.
Considering the differential system \eqref{semi_discrete} and for given degrees of freedom, $\ba^{n}$ represents the vector field which approaches $\varphi$ at time $t^{n}$.
To  increase the accuracy in time, we focus on composition methods for certain numerical methods: the backward Euler and Heun's methods, because of their simplicity and accuracy, designated respectively by the acronyms BE and HM.
\\
On the one hand, we approximate the solution at time $t^{n}$ using the backward Euler scheme as follows:
\begin{equation*}
\varPhi^{\text{BE1}}_{\Delta t}\big(\ba^{n-1}\big) := \ba^{n} =
 \Big(  {\bf{M}}(t^{n}) + \Delta t   {\bf{K}}(t^{n})\Big)^{-1} \,   {\bf{M}}(t^{n-1}) \ba^{n-1}.
\end{equation*}
On the other hand, we use the one-step Heun's method as stated by the following numerical flow:
\begin{equation*}
 \varPhi^{\text{HM1}}_{\Delta t}\left(\ba^{n-1}\right) := \ba^{n} =
 \Bigg(
     \bI - \frac{\Delta t}{2}
              \Big( 
                       \left[ {\bf{M}}\left(t^{n}\right)\right]^{-1}   {\bf{K}}\left(t^{n-1}\right)
                      +  \left[ {\bf{M}}(t^{n})\right]^{-1}   {\bf{K}}\left(t^{n}\right)\left(  \bI - \Delta t \left[ {\bf{M}}(t^{n})\right]^{-1}   {\bf{K}}\left(t^{n-1}\right)
                                                                                                                                \right)
              \Big)
\Bigg)\ba^{n-1}.
\end{equation*}


\begin{figure}[!t]
 \centerline{
     \includegraphics[height=5.2cm]{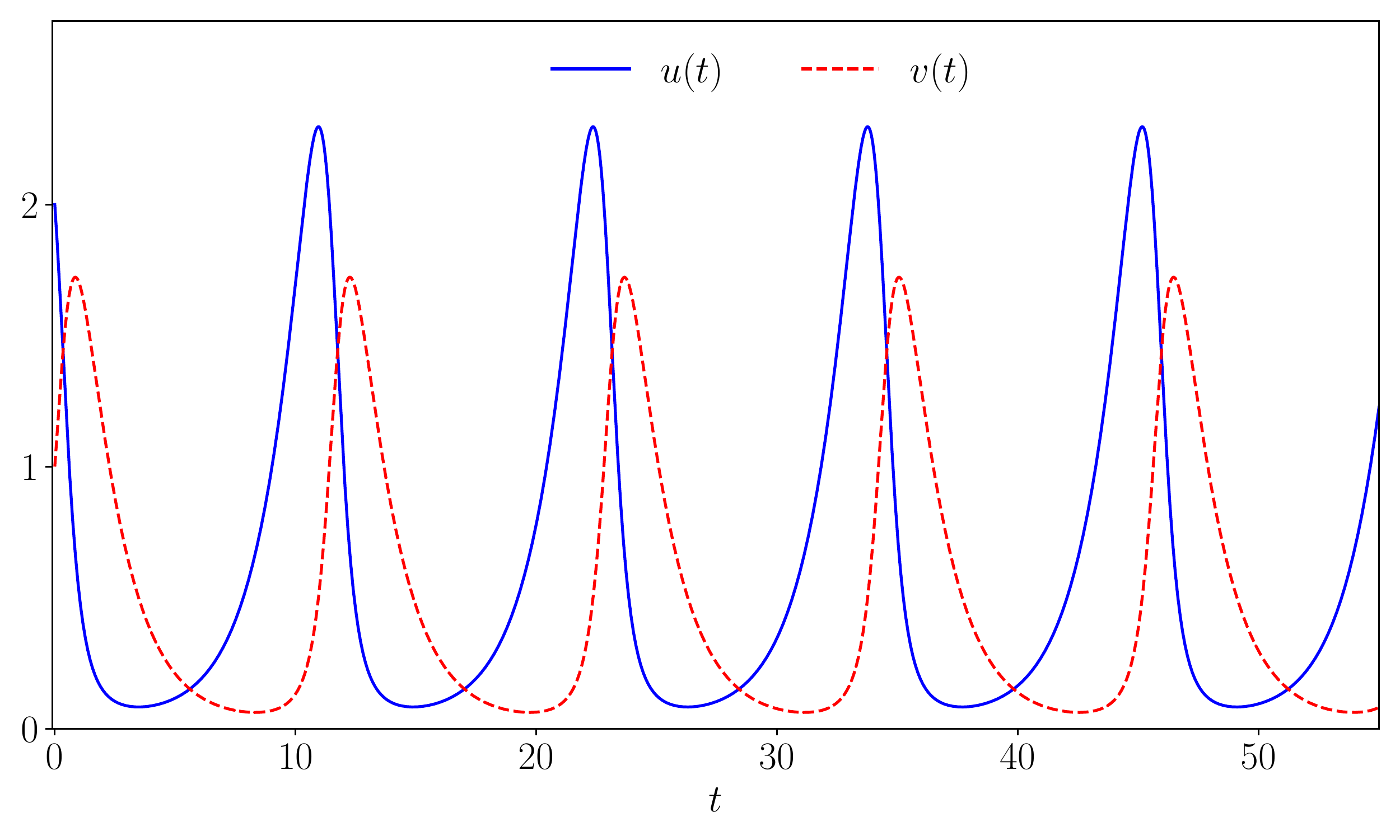}
    \hspace{1.5cm}
     \includegraphics[height=5.9cm]{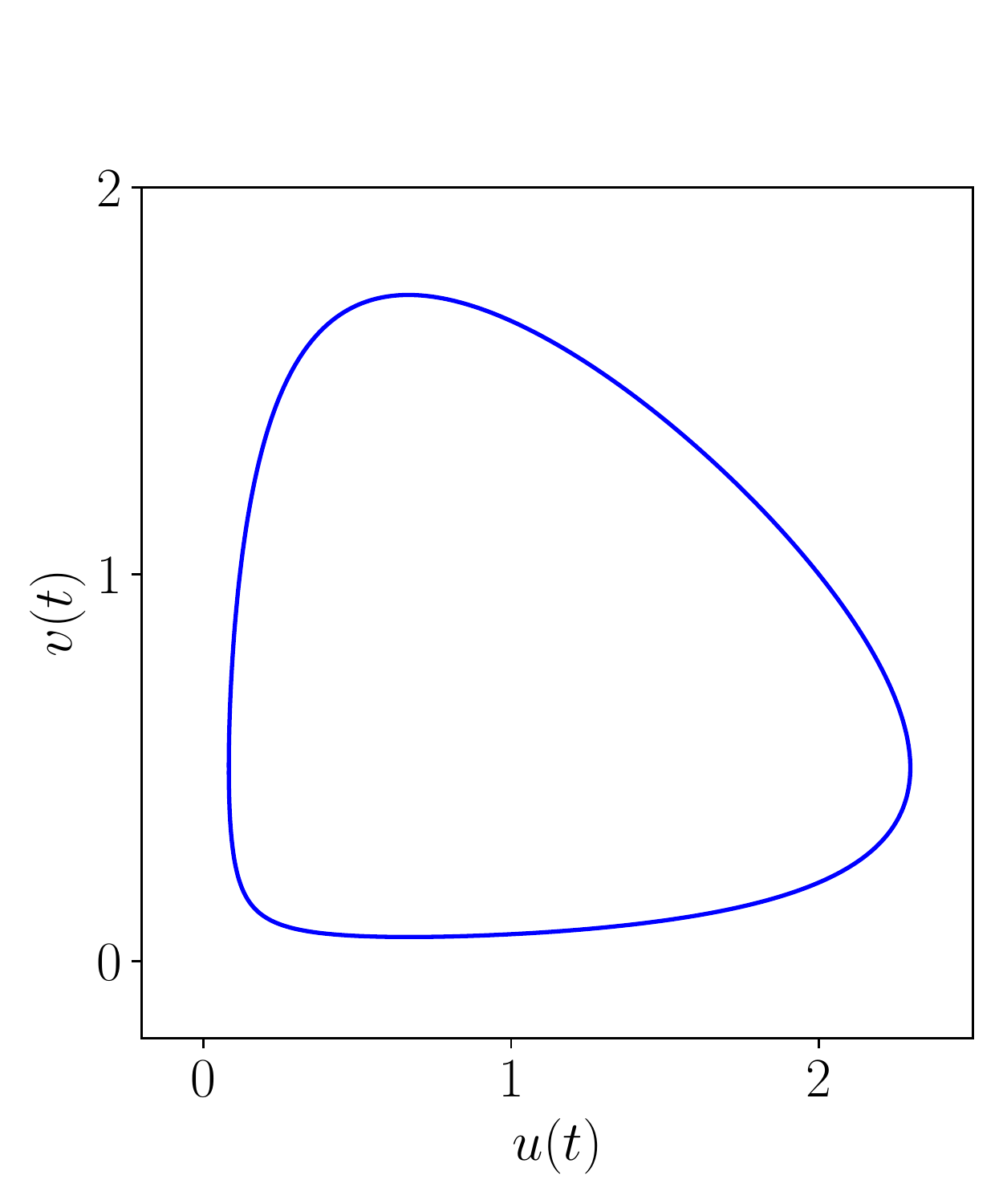}}
 \caption{Numerical solution of Lotka-Volterra equations for the first five periods obtained with the BE4 scheme.}
 \label{fig:equiAngle}
\end{figure}
\section{Numerical examples}\label{sec:sim}
%
\subsection{Example 1: Validation in the case of Lotka–Volterra equations}
\cblue{We apply the above numerical flow compositions to the Lotka-Volterra system of first order nonlinear differential equations describing the evolution over time of two species: predators and prey.} Dynamically speaking, if we assume that predators and prey are not connected, then the number $u$ of prey increases exponentially, with a factor $\alpha$, and $v$ for predators decreases exponentially, with a factor $\delta$.
\cblue{The Lotka-Volterra system describes the predator-prey interactions, when connected, as follows:}
\begin{equation*}
\label{lotka_Volterra}
\left\lbrace
\begin{array}{rcl}\displaystyle
{u}^{\prime}(t) =  &\;\; \alpha ~ u(t) & - \; \beta ~ u(t)v(t)\\
{v}^{\prime}(t) = &- \delta ~ v(t) & + \; \gamma ~ u(t)v(t)\\
\end{array} \right.,
\quad \text{with initial conditions }\quad u(0) = u_0, \, v(0) = v_0.
\end{equation*}
Here, $\beta$ represents the rate of prey mortality by predators and $\gamma$ is the growth rate of predators while eating prey.
%
We present in Fig. \ref{fig:equiAngle} the evolution of the numerical solution of the Cauchy problem $(u(0) = 2, v(0) = 1)$ in the interval $[0,55]$ with $\alpha=\delta=\frac{2}{3}$, $\beta=\frac{4}{3}$ and $\gamma=1$.
\cblue{This problem has an invariant prime integral $F_t$, which is therefore identified using the initial conditions, i.e. $F_t\big(u(t),v(t) \big) = F( u_0,v_0 )$, for any time $t>0$.} 
The invariant is given by:
$$
F_t\big(u(t),v(t)\big) = \beta v(t) + \gamma u(t) - \alpha \log(v(t)) - \delta \log(u(t)).
$$
To study the convergence properties of the numerical strategy, we calculate the error between the approximate invariant and the reference invariant at the initial time $F_{\text{0,ref}}(u_0,v_0)$. The error quantification corresponds to the time evolution of the quantity $F_{\text{t}}(u(t),v(t))$ given by:
\begin{equation}
\label{error_prime_integral}
\|e_{\Delta t^l}\|=\left(\frac{\sum_{n=1}^{N}\big|F_{\text{0,ref}} (u_0,v_0)-F_{\text{t}}(u(t^n),v(t^n))\big|^2 }{\sum_{n=1}^{N}|F_{\text{0,ref}}(u_0,v_0)|^2}\right)^{1/2}
\quad
\mbox{and}\quad
\text{ROC} = \frac{\log_{10}\left(\frac{\|e_{\Delta t^l}\|}{\|e_{\Delta t^{l-1}}\|}\right)}{
\log_{10}\left(\frac{\Delta t^l}{\Delta t^{l-1}}\right)}.
\end{equation}
The parameter ROC represents the convergence rates necessary to set a computational effort to establish a certain precision, in which $\Delta t$ is the time step in a given time discretization level $l$.
To investigate the accuracy of different composition methods, we run a series of simulations using the methods of Backward Euler and Heun and their recursive composition for several time steps.
The table \eqref{tab:p2-t1} and the figure \ref{fig:error-t} report the errors \eqref{error_prime_integral} over a complete period for different compositions of numerical flows.
The estimated convergence orders are reported in the ROC column, 
showing that we retrieve the expected theoretical convergence orders of the basic schemes.
Moreover, the increase in order by double composition of flows is clearly depicted.

\begin{table}[!t]
{
 \begin{tabular}{lcccccccccccc}
 \hline\noalign{\smallskip}
 $\Delta t$& $e^{[BE1]}$ &ROC & $e^{[BE2]}$ &ROC & $e^{[BE3]}$ &ROC & $e^{[HM1]}$ &ROC & $e^{[HM2]}$ &ROC & $e^{[HM4]}$ &ROC \\[1ex]
 \hline\hline\noalign{\smallskip}
4.17E-1  &   4.125  &  $--$  &   2.148E-1  &  $--$  &   3.899E-3  &   $--$  &   2.051E-1  &   $--$  &   6.650E-3  &   $--$  &   2.888E-4  &   $--$ \\[0.5ex]
2.04E-1 & 1.022 & 1.955 & 5.231E-2 & 1.979 & 3.761E-4 & 3.277 & 4.476E-2 & 2.132 & 6.021E-4 & 3.365 & 1.809E-5 & 3.881  \\[1ex]
1.01E-1 & 4.156E-1 & 1.279 & 1.270E-2 & 2.013 & 4.002E-5 & 3.186 & 1.046E-2 & 2.067 & 6.338E-5 & 3.201 & 1.129E-6 & 3.945  \\[1ex]
5.03E-2 & 1.899E-1 & 1.122 & 3.113E-3 & 2.013 & 4.595E-6 & 3.100 & 2.527E-3 & 2.034 & 7.253E-6 & 3.105 & 7.056E-8 & 3.971 \\[1ex]
2.51E-2 & 9.101E-2 & 1.057 & 7.698E-4 & 2.008 & 5.503E-7 & 3.051 & 6.211E-4 & 2.017 & 8.667E-7 & 3.054 & 4.412E-9 & 3.985 \\[1ex]
1.25E-2 & 4.457E-2 & 1.028 & 1.914E-4 & 2.005 & 6.733E-8 & 3.025 & 1.540E-4 & 2.009 & 1.059E-7 & 3.027 & 2.758E-10 & 3.992 \\[1ex]
6.25E-3 & 2.206E-2 & 1.014 & 4.770E-5 & 2.002 & 8.328E-9 & 3.013 & 3.833E-5 & 2.004 & 1.309E-8 & 3.014 & 1.725E-11 & 3.996 \\[1ex]
3.13E-3 & 1.097E- & 1.007 & 1.191E-5 & 2.001 & 1.035E-9 & 3.006 & 9.561E-6 & 2.002 & 1.627E-9 & 3.007 & 1.084E-12 & 3.999 \\[1ex]
1.56E-3 & 5.473E-3 & 1.003 & 2.975E-6 & 2.001 & 1.291E-10 & 3.003 & 2.388E-6 & 2.001 & 2.027E-10 & 3.003 & 6.774E-14 & 3.999 \\[1ex]
7.81E-4 & 2.733E-3 & 1.002 & 7.434E-7 & 2.000 & 1.612E-11 & 3.001 & 5.966E-7 & 2.001 & 2.531E-11 & 3.002 & 1.385E-14 &  $--$ \\[1ex]
3.91E-4 & 1.36E-3 & 1.001 & 1.858E-7 & 2.000 & 2.013E-12 & 3.002 & 1.491E-7 & 2.000 & 3.156E-12 & 3.003 & 2.510E-14 &  $--$ \\[1ex]
1.95E-4 & 6.827E-4 & 1.000 & 4.645E-8 & 2.000 & 2.231E-13 & 3.173 & 3.727E-8 & 2.000 & 3.710E-13 & 3.088 & 5.621E-14 & $--$ \\[1ex]
9.73E-5 & 3.400E-4 & 1.000 & 1.152E-8 & 2.000 & 6.210E-14 & 1.835 & 9.245E-09 & 2.000 & 4.395E-14 & 3.060 & 8.845E-15 & $--$ \\[1ex]
 \multicolumn{2}{>{\columncolor{cgray!0}}l}{Order of convergence}& 1 & $--$ &  2 & $--$ &  3& $--$&2& $--$&3& $--$&4 \\[0ex]
 \hline\hline
 \end{tabular}
}
 \caption{Time convergence history for different time integration composition methods.}
 \label{tab:p2-t1}
\end{table}

\begin{figure}[!h]
 \centerline{
     \includegraphics[height=5.2cm]{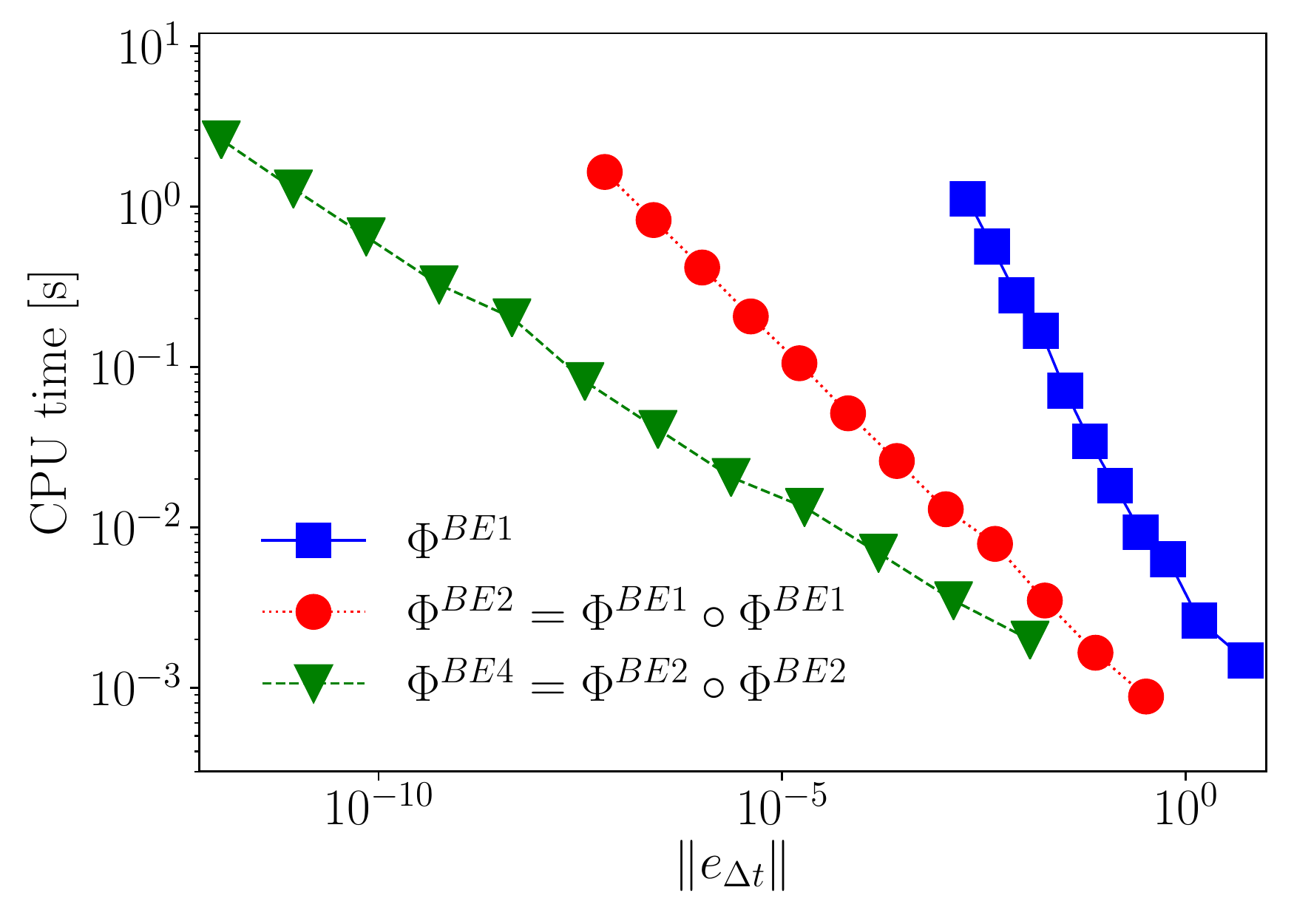}
     \hspace{0.75cm}
     \includegraphics[height=5.2cm]{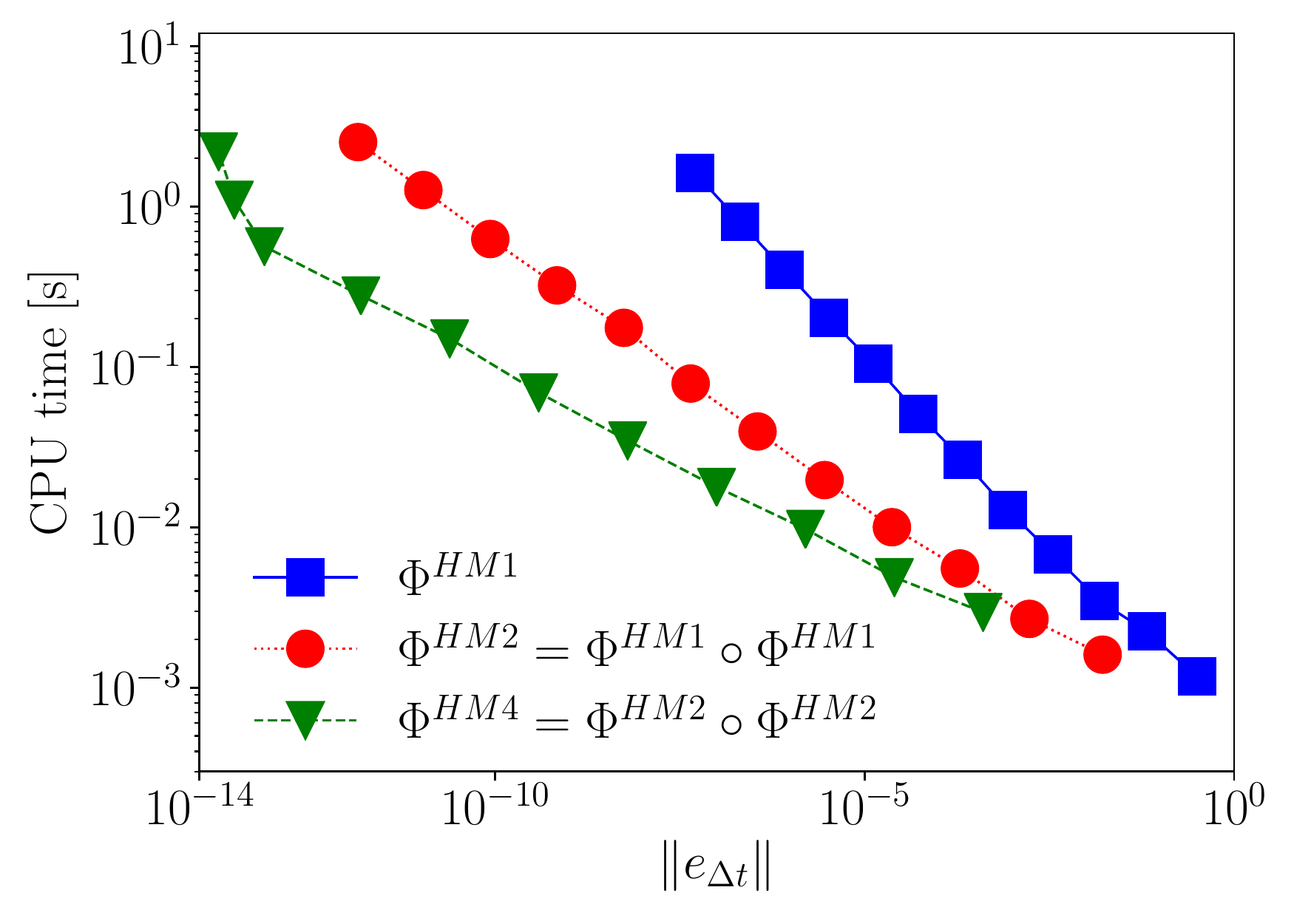}}
 \caption{CPU time versus error for different compositions of backward Euler (left) and Heun (right) schemes.}
 \label{fig:CPU-Error}
\end{figure}

Thereafter, we study the computational efficiency by evaluating the CPU times for different time steps and composition methods. Thus, a correlation can be established between the error sought by the user and the CPU time that the numerical flow needs to carry out the simulation, for the two discrete schemes BE and HM and their compositions.
Results in Fig. \ref{fig:CPU-Error} show that \cred{the resulting composition scheme allows the user to produce numerical results with a lower computational time than the original scheme, while keeping the same precision. This has been observed for different composition levels, for a given original basic scheme (Heun or Backward Euler method).}
\cblue{Note that we do not compare in this work such compositions with other schemes having the same convergence orders.}
%
%

%

\subsection{Example 2: Level set problem by composition methods}


\cblue{In this example, the aforementioned composition of one-step time integration schemes is used for simulations of largely deformable interfaces in a level set framework.}
The numerical computations are performed using the open-source software FEniCSx\cite{AlnaesEtal2015}.
We first consider a reversible vortex test case and evaluate the properties of convergence with respect to the temporal discretization.
The computational domain is $\Lambda=[0,1]^2$.
An initially circular interface of radius $R=0.15$ is centered at $(0.7,0.7)$ and is periodically stretched into thin filaments by a vortex flow field.
\cblue{The interface unrolls and reaches its maximum deformations at $t=T/2$, before resuming its circular shape at time $T=4$}.
The advection velocity \cblue{is given by}:
\[
\bu(t,\bx) =\left( 
-2 \sin( \pi x)^2  \ sin(\pi y) \cos( \pi y) \cos(\pi t/T)
\,,\,
 2 \sin(\pi y)^2   \sin( \pi x) \cos( \pi x) \cos(\pi t/T)
 \right)^T,
 \quad\text{with}\,\,
 \bx=(x,y)^T\in\Lambda.
\]

\begin{figure}[!h]
    \centering
    \begin{minipage}{0.68\textwidth}
     \hspace{0.0cm}
     \includegraphics[width=0.23\textwidth]{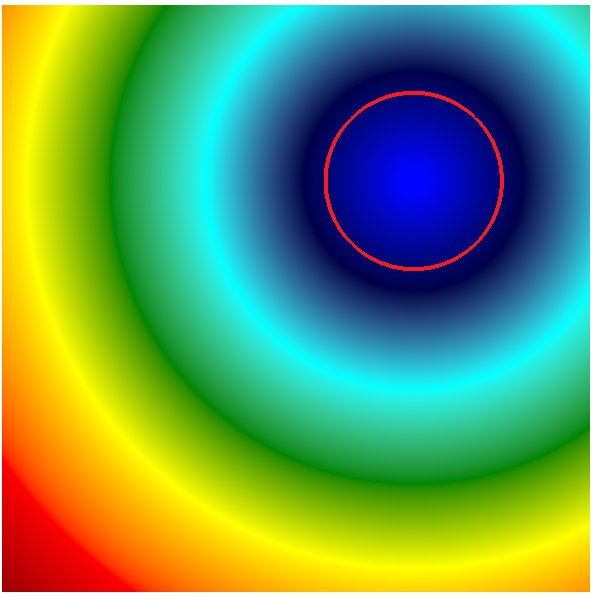}
     \includegraphics[width=0.23\textwidth]{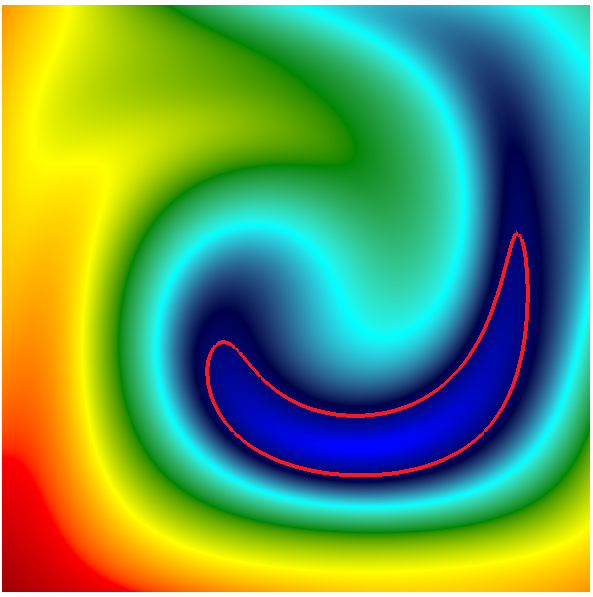}
     \includegraphics[width=0.23\textwidth]{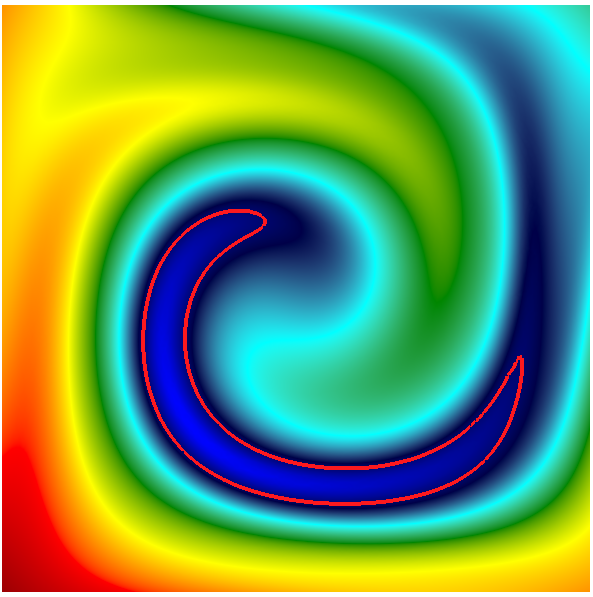}
     \includegraphics[width=0.23\textwidth]{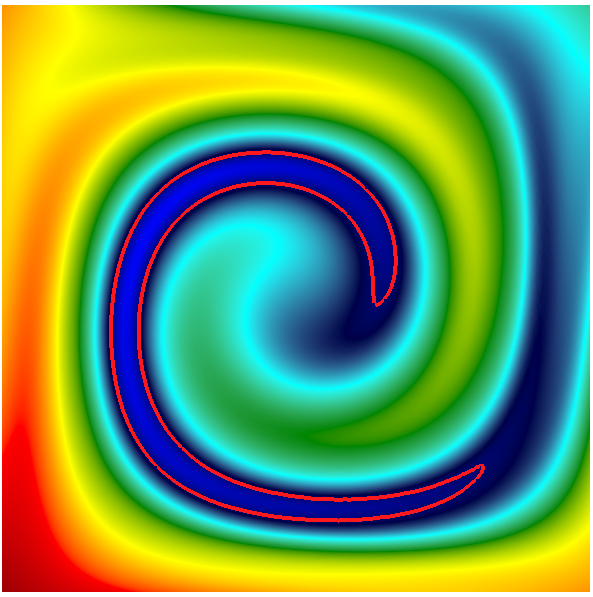}
     \newline
     \includegraphics[width=0.23\textwidth]{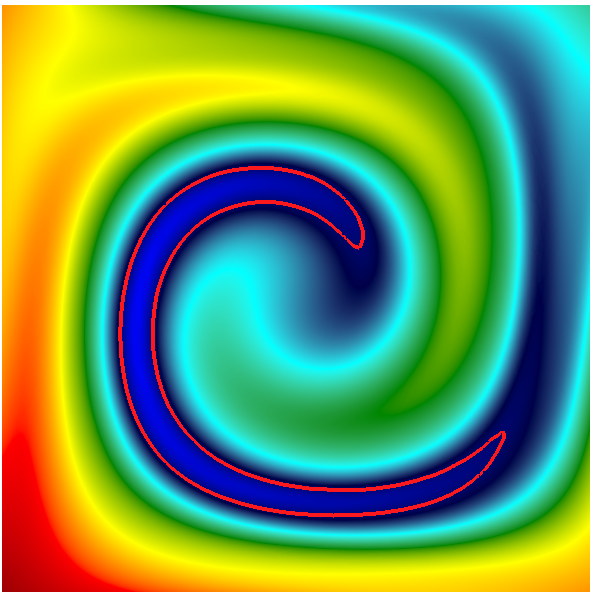}
     \includegraphics[width=0.23\textwidth]{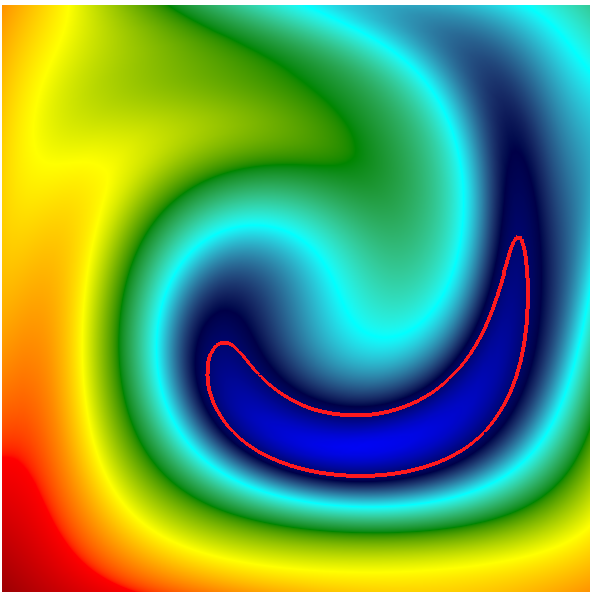}
     \includegraphics[width=0.23\textwidth]{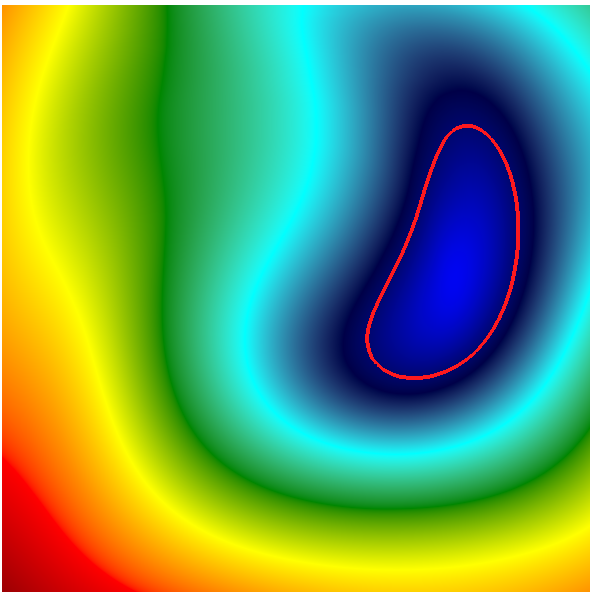}
     \includegraphics[width=0.23\textwidth]{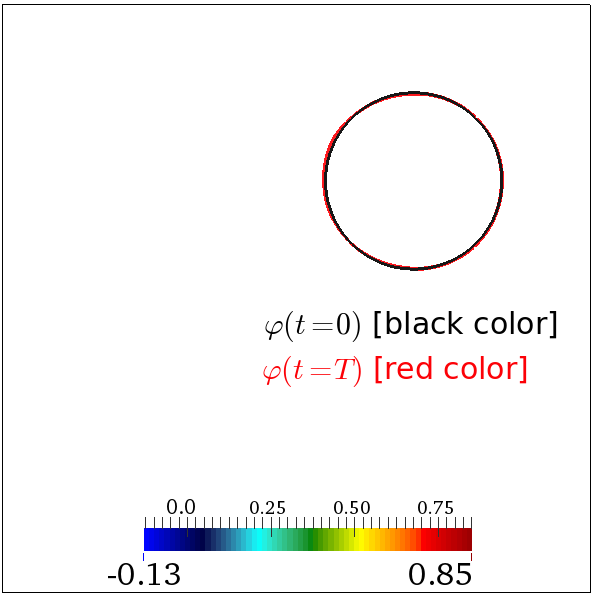}
    \end{minipage}
    \begin{minipage}{0.29\textwidth}
    \hspace{-0.15cm}
     \includegraphics[height=5.25cm]{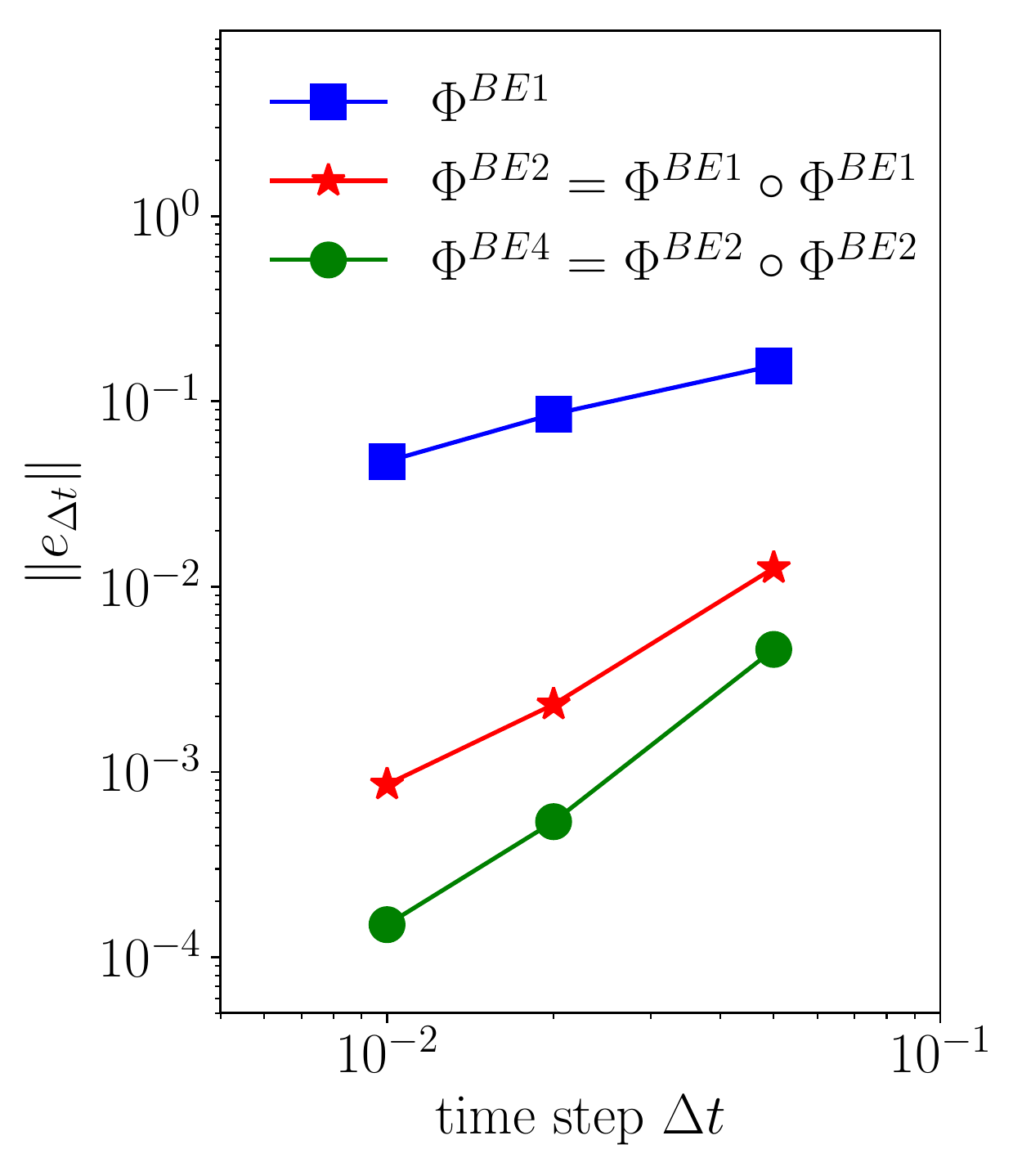}
    \end{minipage}
\caption{(Left) Snapshots showing interface deformations at 
$t\in\{0, 0.6, 1, 2,2.6, 3.4, 3.8, 4\}$, for $h=10^{-2}$.
(Right) Convergence history for various time steps and time composition schemes.}
\label{fig:test}
\end{figure}

Some snapshots showing the interface deformations are reported in Fig. \ref{fig:test}.
\cblue{We investigate the order of accuracy with respect to the time discretization using a relatively simple composition method. 
Let $\pi_h$ be the Lagrange interpolation operator.
We compute the errors in the $L^2(\Lambda)$ norm for different time step sizes with respect to an exact reference solution $\pi_h\varphi$ at time $t=T$ after one deformation period.
The convergence of the composition technique is assessed using the implicit backward Euler scheme, as well as two- and four-time backward Euler compositions, respectively called $\varPhi^{\text{BE2}}$ and $\varPhi^{\text{BE4}}$.
Fig \ref{fig:test} (right) reports the convergence of the computed errors versus the time discretization, showing in particular an increase in the order of convergence by composition of discrete flows.
}

\begin{figure}[!h]
    \begin{minipage}{0.63\textwidth}     
    \centerline{
     \includegraphics[width=0.24\textwidth]{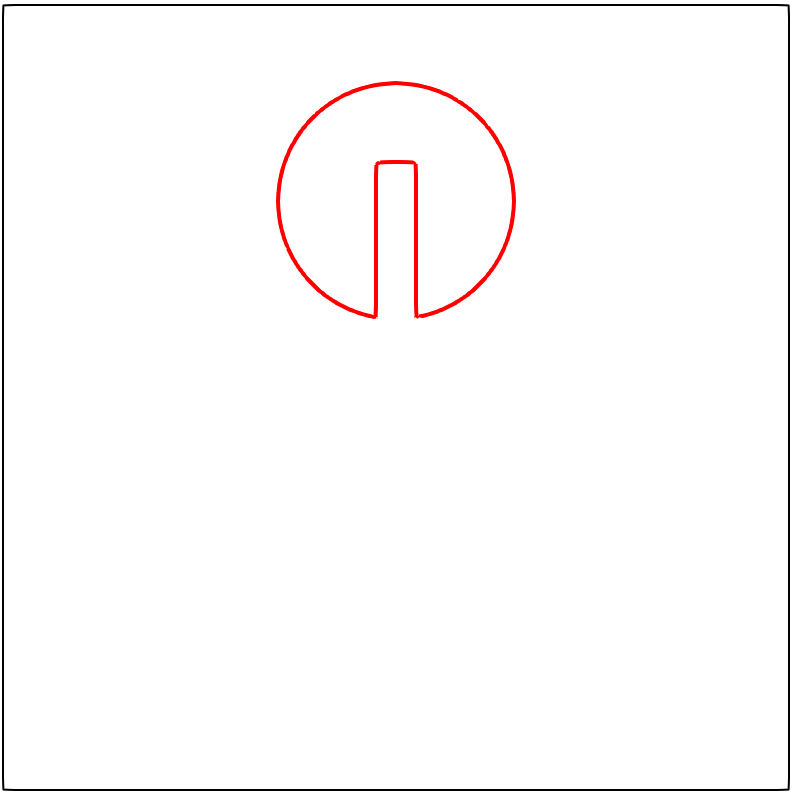}
     \includegraphics[width=0.24\textwidth]{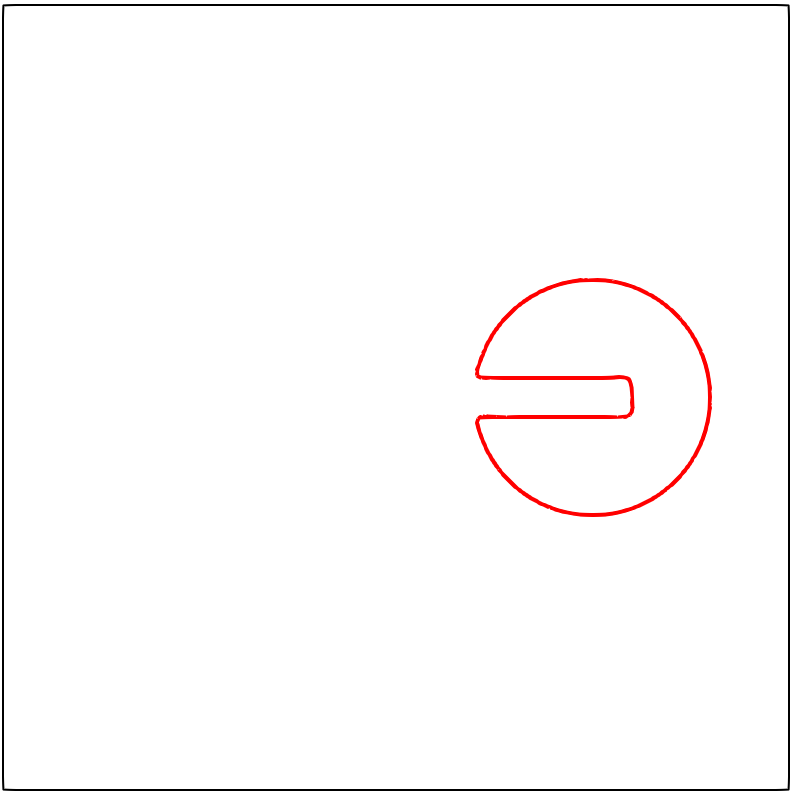}
     \includegraphics[width=0.24\textwidth]{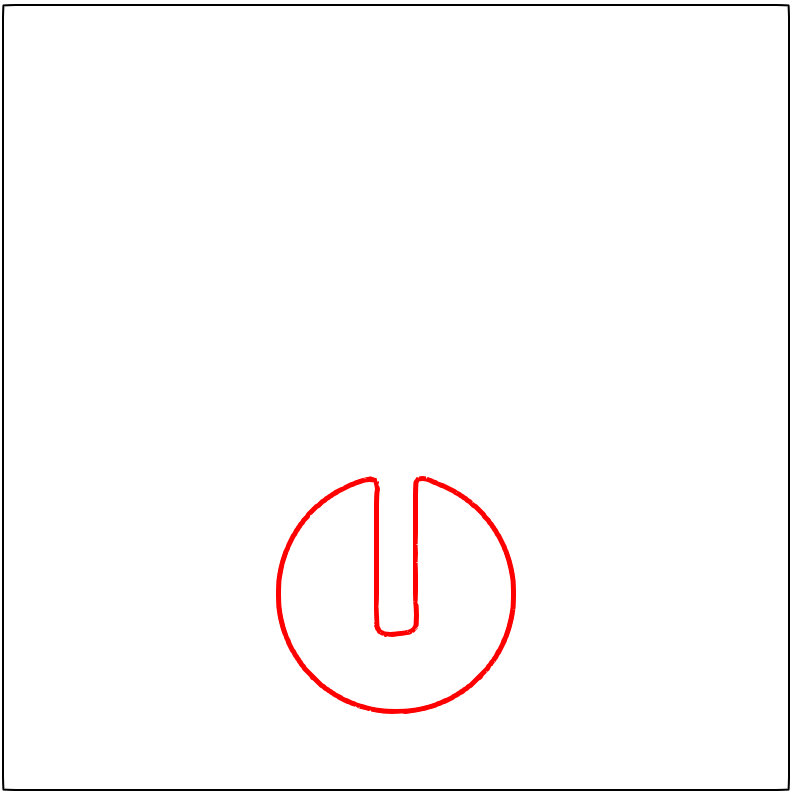}
     }
     \centerline{
     \includegraphics[width=0.24\textwidth]{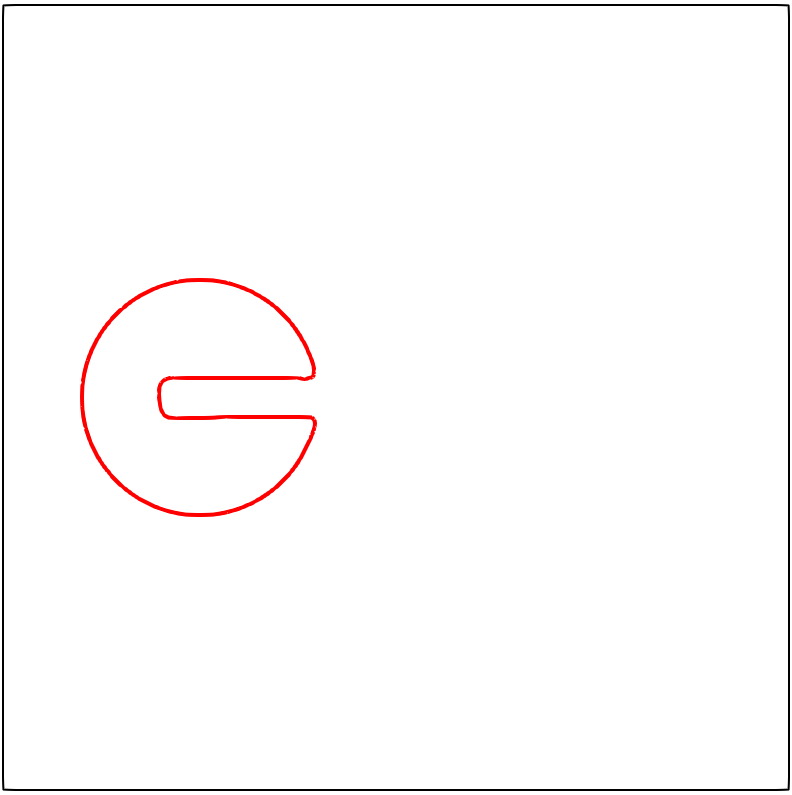}
     \includegraphics[width=0.24\textwidth]{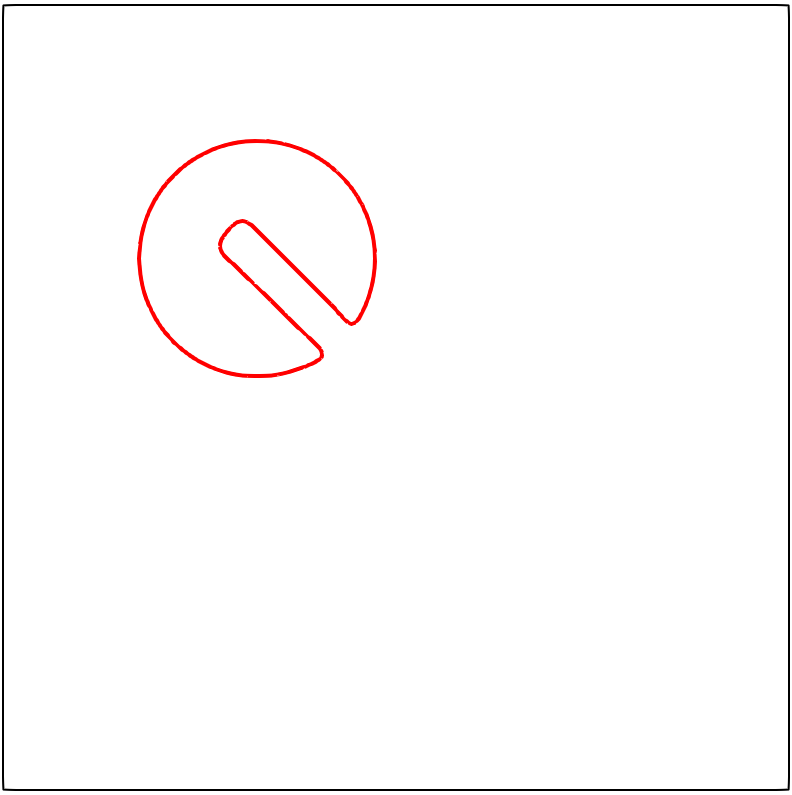}
     \includegraphics[width=0.24\textwidth]{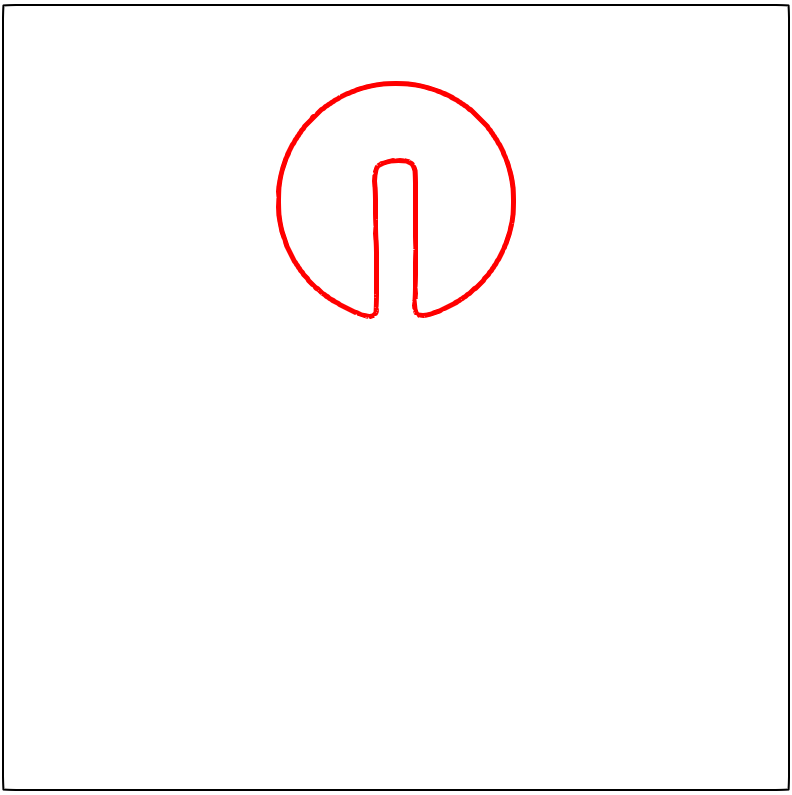}
     }
    \end{minipage}
    \begin{minipage}{0.35\textwidth}     
     \includegraphics[width=0.86\textwidth]{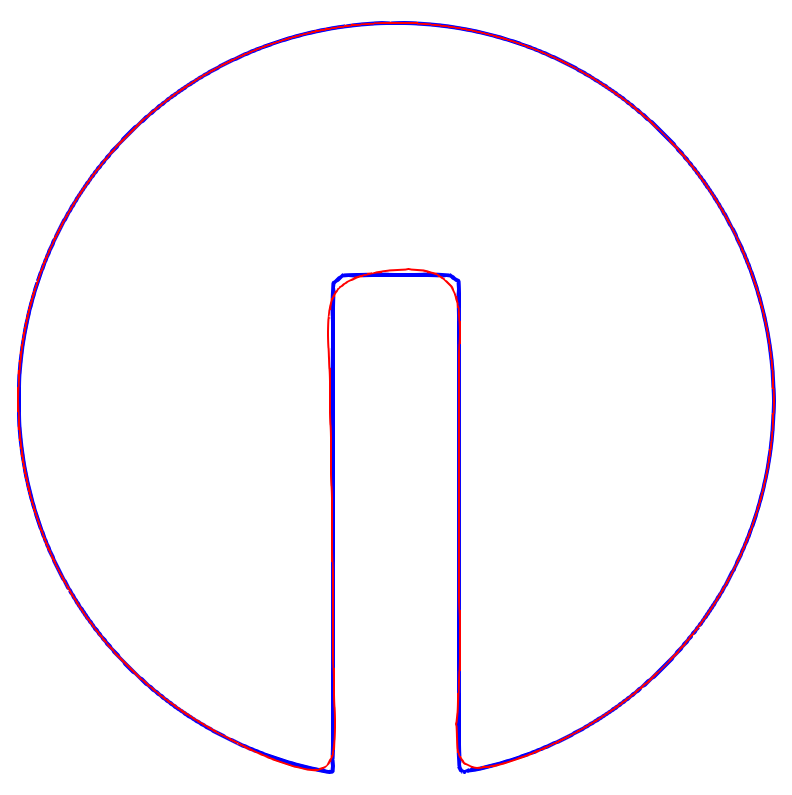}
    \end{minipage}
     \caption{\cblue{Zalesak's rotating disk. Snapshots showing interface deformations at 
$t\in\{0, 0.6, 1, 2,2.6, 3.4, 3.8, 4\}$.}} 
\label{fig:zal}
\end{figure}

Finally, we assess the robustness of our level set solver in the case of the standard Zalesak's rotating disk test, using a one-time composition of the backward Euler scheme.
We consider the same setting presented in \cite{sar16}, where the slotted disk recovers the initial position after a period $T=4$. We choose $\Delta t=10^{-3}$ and $h=10^{-2}$.
A few snapshots of the rotating disk are shown in Fig \ref{fig:zal}, showing good preservation of the sharp angles.

\vspace{-0.2cm}
\section{Conclusion}
We have presented an application of numerical methods of flow composition to the simulation of the level set problem in order to increase the accuracy of the solution and to capture the dynamics of stiff problems.
A validation is presented in the case of Lotka-Volterra differential equations.
By using simple and lower-order numerical schemes, we have shown numerically that composition techniques allow to increase the accuracy of the approximation while presenting a lower computational cost.
This is part of an ongoing work on modeling the dynamics of highly-deformable biomembranes \cite{LSMS2018,Gizzi2015157,Kolahdouz20157}, where temporal integration schemes with high orders are required.
It is also planned to implement the method for fluid-structure interaction problems with large structural deformations \cite{LS_JCAM}.



\begin{acknowledgments}
The authors acknowledge financial support from KUST through the grant FSU-2021-027.
\end{acknowledgments}

\bibliography{mybibfileAML}

\end{document}